\documentclass[11pt]{article}

\usepackage{graphicx}

\usepackage[pdfdisplaydoctitle=true]{hyperref}

\usepackage{xcolor}

\usepackage{amssymb}
\usepackage{framed}
\usepackage{fullpage}
\usepackage{color}
\usepackage{setspace}
\usepackage{amsmath}
\usepackage{helvet}

\usepackage[letterspace=20]{microtype}

\SetTracking[ spacing = {2500000*,166, } ]{ encoding = *, shape = sc }{12000000}

\usepackage{titlesec}
\titleformat*{\section}{\fontsize{18}{22}\selectfont\bfseries}
\titleformat*{\subsection}{\fontsize{16}{20}\selectfont\bfseries}
\titleformat*{\subsubsection}{\fontsize{14}{17}\selectfont\bfseries}

\hypersetup{
    pdftitle={Asymptotics of two--point correlations in the multi--species $q$--TAZRP}, 
    pdfauthor={Jeffrey Kuan and Zhengye Zhou}, 
    pdfsubject={Research Paper},  
    pdfkeywords={Multi-species, correlations, q-TAZRP},  
    pdfcreator={Jeffrey Kuan},  
    bookmarksnumbered=true,     
    bookmarksopen=true, 
    linktocpage=true, 
} 

\usepackage[utf8]{inputenc}
\usepackage{amsmath,amssymb,amsfonts,amsthm}
\usepackage{graphicx}
\usepackage{fullpage}
\usepackage{xcolor}
\newtheorem{thm}{Theorem}

\theoremstyle{remark}
\newtheorem{remark}[thm]{Remark}

\usepackage{fullpage}

\begin{document}

\title{Asymptotics of two--point correlations in the multi--species $q$--TAZRP}
\author{Jeffrey Kuan and Zhengye Zhou}
\date{}

\Large

\maketitle

\large 

\begin{center}

\fbox{\parbox{0.8\textwidth}{Accessibility Statement: \\~\\  A WCAG2.1AA compliant version of this PDF is available  at \href{http://math.tamu.edu/~jkuan/arXiv-2305.00321.html}{the first author's professional website.} The accessible version of the PDF was created using LuaLaTeX; arXiv does not support LuaLaTeX nor does it allow for PDFs generated by LuaLaTeX to be directly uploaded.  \\~\\  The accessible version of this PDF has alternative text describing the equations and figures, which may be useful to the reader.\\~\\ This PDF was created in consultation with the British Dyslexia Association's Style Guide 2023, although arXiv moderators did not allow me to upload the version with better color contrast.  \\~\\  For any issues accessing the information in this document, please \href{mailto:jkuan@math.tamu.edu}{email the first author.}}}

\end{center}

\abstract{\large  A previous paper by the authors found explicit contour integral formulas for certain joint moments of the multi--species $q$--TAZRP (totally asymmetric zero range process), using algebraic methods. These contour integral formulas have a ``pseudo--factorized'' form which makes asymptotic analysis simpler. In this brief note, we use those contour integral formulas to find the asymptotics of the two--point correlations. As expected the term arising from the ``shift--invariance'' makes a non--trivial asymptotical contribution.}

\large 

\section{Introduction}
In interacting particle systems, control of the two-point correlations are generally needed to determine the hydrostatic limit. For symmetric exclusion processes, this has been done in \cite{BMNS} and \cite{Kyoto}. For symmetric inclusion process, this was done in \cite{FGS22}. In asymmetric processes, calculating the correlations is naturally harder, such as \cite{UW05}. 

In more recent years, a certain type of interacting particle systems called ``multi--species'' models have been well--studied, such as \cite{BS,BS2,BS3,BGW,Gal,Kuan-IMRN,KuanCMP,KMMO,Take15,ZZ}. In this paper, we consider one specific asymmetric multi--species model, called the multi--species $q$--TAZRP. In previous work \cite{FKZ}, the authors and C. Franceschini found an explicit expression for the two--point correlations of the height function of this model. (That work originated with orthogonal polynomial duality (e.g. \cite{ACR18,ACR21,CFG21,CFGGR,FG,Gro19,ZZ}), but this current paper will not reference duality). This expression is written as the sum of six terms, with each term being a product of a simple $q$--Pochhammer (which only depends on the number of initial particles) and a contour integral. Of these six terms, five can be derived from the Markov projection (or color--blind) property, but the sixth used a so--called ``shift--invariance''. 

The contour integral and product form of each term makes the asymptotics suitable for analysis, using standard and well--known methods. In this brief note, we carry out these calculations and make the observation that the term arising from shift--invariance makes a non-trivial asymptotic contribution, as expected.

Acknowledgements. The authors would like to thank Chiara Franceschini for helpful discussions at the workshop ``Recent Developments in Stochastic Duality'' in Eindhoven, and for support from EURANDOM and the NWO grant 613.001.753 ``Duality for interacting particle systems''. The first author was supported by National Science Foundation grant DMS-2000331.

\section{Main Result}

\subsection{Definition of Model}
The multi--species $q$--TAZRP will be defined informally here. It was first introduced by Takeyama, generalizing the single--species version by Sasamoto and Wadati. There are $n$ different species of particles, labeled $1,2,\ldots n$. Arbitrary many particles of any number of species may occupy a lattice site. The lattice is the one--dimensional lattice, labeled $\mathbf{Z}$, and particles jump one step to the right at a time. 

The jump rates depend only on the number of particles at a single site, as is the case for any zero range process. Let $[k]_q=1-q^k$ denote the $q$--integer in the below equations. If the number of species $k$ particles at a site is denoted by $z_k$, then the jump rate for the species $k$ particle is $q^{z_1 + \cdots + z_{k-1}} [z_k]_q{/(1-q)}.$ 

\begin{figure}[h]
\includegraphics[height=3in]{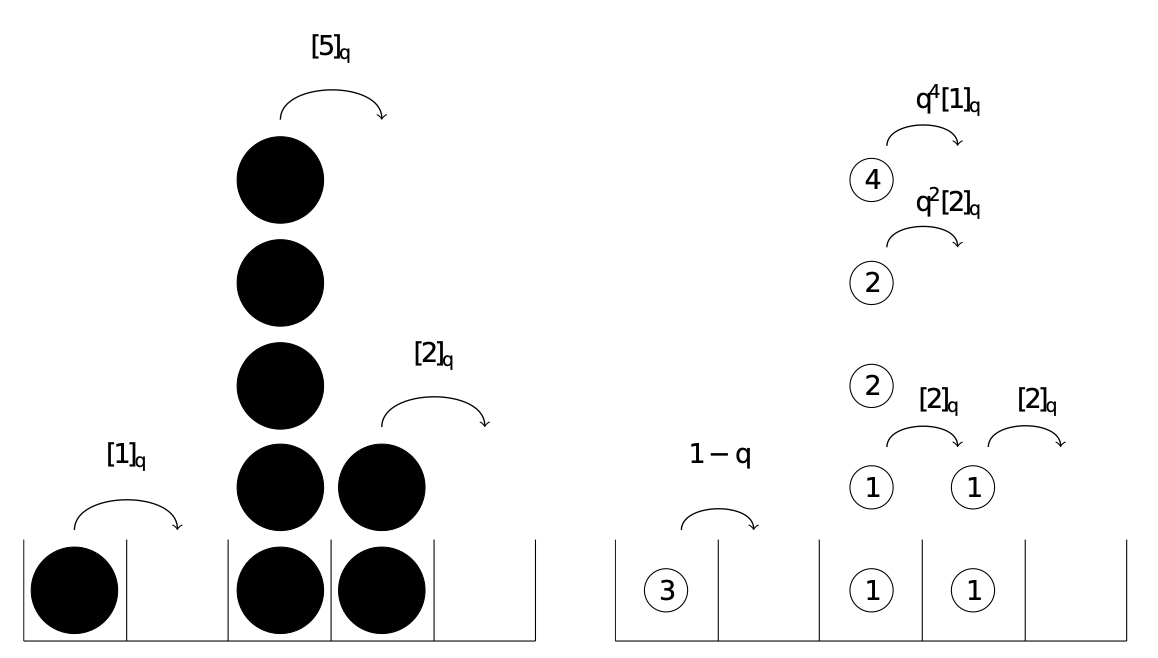}
\caption{The jump rates for the $q$--TAZRP.}
\label{image}
\end{figure}

Based on these jump rates, the multi--species $q$--TAZRP has a ``Markov projection'' property: if one only ignores the species completely, one obtains the (single--species) $q$--TAZRP. See the Figure \ref{image} for an example of the jump rates. 

Recall the incomplete gamma function $$\Gamma(a,z) =\int_{z}^{\infty} t^{a-1} \mathrm{e}^{-t} \mathrm{~d} t
$$
and define its normalization 
$$
Q(a, z)=\frac{\Gamma(a, z)}{\Gamma(a)} .
$$
Recall that for a continuous--time totally asymmetric simple random walk $X(t)$ starting at $0$, with jump rates $r$, then 
$$
\mathbf{P}(X(t) \leq n) = \sum_{k=0}^n e^{-rt}\frac{(rt)^k}{k!} = \frac{\Gamma(n+1,rt)}{n!} = Q(n+1,rt).
$$

We also have
$$ 
{ }_{1} \phi_{0}\left(\begin{array}{c}q^{-n} \\ -\end{array} ;, z\right)=\left(z q^{-n} ; q\right)_{n}, \quad n=0,1,2, \ldots
$$
where the $q$--Pochhammer is
\begin{equation*}
(a ; q)_{n}=(1-a)(1-a q) \cdots\left(1-a q^{n-1}\right), \quad(a ; q)_{0}:=1.
\end{equation*}

In Theorem 3.5 of \cite{FKZ}, it was proved that

\begin{thm}\label{Application}
Suppose $\boldsymbol\xi$ has initial conditions with $n_1$ species 0 particles at site $x_1$ and $n_2$ species 1  particles at site $x_2$, or in other words 
$$
\boldsymbol\xi_i^x 
= 
\begin{cases}
n_1 \ \text{for} \  i=0,x=x_1,\\
n_2 \ \text{for} \  i=1,x=x_2,\\
0 \ \  \text{else},
\end{cases}
$$  
and $\boldsymbol\xi$ evolves as a $q$--TAZRP with total asymmetry to the right. Let $y_1$ and $y_2$ be some lattice sites, and additionally assuming that $x_1<x_2$. 

The dual process $\boldsymbol\eta$ has 1 species 0 particles at site
$y_1$, 1 species 1 particles at site $y_2$, so in other words
$$
\boldsymbol\eta_i^y 
= 
\begin{cases}
1 \ \text{for} \ i=0,y=y_1,\\
1 \ \text{for} \ i=1,y=y_2,\\
0 \ \  \text{else}.
\end{cases}
$$  
Then there is the explicit formula 
$$
\mathbb{E}_{\boldsymbol{\xi}}[ q^{1/2h(\boldsymbol\xi (t),\boldsymbol\eta)} \prod_{i=0}^{1}\prod_{x=1 }^{L}\ _1\phi_0\left(q^{-\xi_i^x (t)};q,q^{-\left(N_{x-1}^-(\boldsymbol\xi_i (t))+N_{x+1}^+(T(\boldsymbol\eta)_{i+1})\right)+1/2}\right) ]  = \sum_{i=1}^6{\delta_ip_i},
$$
where 
$$
h(\boldsymbol{\xi}, \boldsymbol{\eta})=\sum_{x \in \Lambda_{L}} \sum_{i=0}^{n-1}\left(-\xi_{i}^{x} N_{x+1}^{+}\left(\boldsymbol{\eta}_{[0, n-2-i]}\right)+\eta_{i}^{x} N_{x}^{+}\left(\boldsymbol{\xi}_{[0, n-2-i]}\right)\right)
$$
\begin{align*}   
   \delta_1:= q^{n_1/2}\left(q^{-n_1-1/2};q\right)_{n_1} \left(q^{-n_2+1/2};q\right)_{n_2}    ,\\
    \delta_2:=     q^{n_1/2}\left(q^{-n_1+1/2};q\right)_{n_1} \left(q^{-n_2+1/2};q\right)_{n_2}  ,\\
     \delta_3:=   q^{-n_1/2}\left(q^{-n_1-1/2};q\right)_{n_1} \left(q^{-n_2+1/2};q\right)_{n_2}    ,\\
      \delta_4:=    q^{-n_1/2}\left(q^{-n_1+1/2};q\right)_{n_1} \left(q^{-n_2+1/2};q\right)_{n_2}  ,\\
   \delta_5:=   q^{-n_1/2}\left(q^{-n_1-1/2};q\right)_{n_1} \left(q^{-n_2-1/2};q\right)_{n_2}  ,\\
    \delta_6:=     q^{-n_1/2}\left(q^{-n_1+1/2};q\right)_{n_1} \left(q^{-n_2-1/2};q\right)_{n_2}  ,
  \end{align*}
  and
  $$
p_1 = q_1-q_4, \quad p_2=q_4, \quad p_3 = q_2+q_3-q_5-q_6, \quad p_4 = -q_3 + q_5+q_6, \quad p_5=q_5, \quad p_6=q_3-q_5,
$$
with
$$
q_1 ={\delta_{y_1\le x_1+}\delta_{y_1> x_1} }Q(y_1-x_1,t) , \quad q_3 = 1 - {\delta_{y_1\ge x_2}Q(y_1-x_2,t)}{-\delta_{y_1<x_2}} , \quad q_2 = 1- q_1-q_3
$$
and
\begin{align*}
q_4&=\frac{ 1}{(2\pi i)^2}\int_{C} \frac{d w_{1}}{w_{1}}  \int_{C} \frac{d w_{2}}{w_{2}} \frac{w_1-w_2}{w_1-qw_2} \prod_{j=1}^{2}\left[(1-w_j)^{-({y_j-x_1}+1)}e^{-w_{j} t}\right] ,\\
q_5&= \frac{ q}{(2\pi i)^2}\int_{\tilde{C}_{1}}  \frac{dw_1}{w_1}  \int_{\tilde{C}_{2}}  \frac{dw_2}{w_2} \frac{w_1-w_2}{w_1-qw_2}(1-w_1)^{-({y_1-x_2}+1)}(1-w_2)^{-({y_2-x_1}+1)} e^{-(w_1+w_2)t  }   ,\\
q_6 &= 1- q_1 - q\frac{ 1}{(2\pi i)^2}\int_{\tilde{C}_{1}}  \frac{dw_1}{w_1}  \int_{\tilde{C}_{2}}  \frac{dw_2}{w_2} \frac{w_1-w_2}{w_1-qw_2}\prod_{j=1}^2\left[ (1-w_j)^{-({y_j-x_1}+1)} e^{-w_jt  } \right],
\end{align*}
{Here, the $C$ are ``large contours'' that are centered at $0$ and contain $1$, and the $\tilde{C}$ are ``small contours'' where $\tilde{C_2}$ contains $1$ and not $0$, whereas $\tilde{C_1}$ contains $q\tilde{C_2},1$ and not $0$.}
\end{thm}

\begin{remark}
Of these six terms, five of them (except $q_5$) can be derived from the two Markov projections to single--species $q$--TAZRP. The term $q_5$ used the ``shift--invariance'' property of \cite{BGW,Gal,Kua21}, and see also \cite{He}; see the fifth bullet point of section 4.9 of \cite{FKZ}. Note that the $\delta_j$ terms only depend on $n_1,n_2$; while the $q_j$ terms do not depend on $n_1$ or $n_2$ at all. This makes asymptotic analysis significantly easier.
\end{remark}

\begin{remark}
There was a typo in Theorem 3.5 of \cite{FKZ} (arXiv version), which has been corrected here, and will be corrected when the refereed version of the paper is updated on arXiv. The typo occured due to conflicting notation when citing the results of \cite{Kua21} in Lemmas 4.13 and 4.14 of \cite{FKZ} . The correction has been confirmed by all three authors of \cite{FKZ}. 
\end{remark}

\subsection{Main Theorem}
The main theorem of this paper is essentially an asymptotic analysis of the contour integrals. Each term will be analyzed individually. Since $\delta_j$ only depend on $n_j$, those terms do \textbf{not} need to be analyzed.

\begin{thm}
Consider the limit as ${y_j-x_i}= L+c_{ij}L^{1/2} +O(1)$ and $t=L$. Then as $L\rightarrow \infty$,
$$
\begin{aligned}
q_1 &\approx   \frac{1}{2} \operatorname{erfc}\left(- \frac{c_{11}}{2^{\frac{1}{2}}}\right)+\frac{1}{\sqrt{2 \pi a_1}} \exp \left(-\frac{c_{11}^{2}}{2}\right) \sum_{n=0}^{\infty} \frac{C_{n}\left(-c_{11}\right)}{a_1^{n / 2}},\\
q_3 &\approx  1- \frac{1}{2} \operatorname{erfc}\left(- \frac{c_{21}}{2^{\frac{1}{2}}}\right)-\frac{1}{\sqrt{2 \pi a_2}} \exp \left(-\frac{c_{21}^{2}}{2}\right) \sum_{n=0}^{\infty} \frac{C_{n}\left(-c_{21}\right)}{a_2^{n / 2}},
\end{aligned}
$$
where $\mathrm{erfc}$ denotes the complementary error function and $C_n(\eta)$ are polynomials in $\eta$ of degree $3n+2$ and satisfy 
$$
C_0(\eta) = \frac{1}{3}\eta^2 - \frac{1}{3},
$$
and
$$
C_{n}(\eta)+\eta C_{n}^{\prime}(\eta)-C_{n}^{\prime \prime}(\eta)=\eta\left(\eta^{2}-2\right) C_{n-1}(\eta)-\left(2 \eta^{2}-1\right) C_{n-1}^{\prime}(\eta)+\eta C_{n-1}^{\prime \prime}(\eta),
$$
and 
$$
a_i = L+\frac{c_{i1}^2}{2}+c_{i1}\sqrt{L+\frac{c_{i1}^2}{4}}.    
$$

Furthermore,
$$
\begin{aligned}
q_4 & \approx\frac{1 }{2\pi }\int_{-\infty}^{c_{1j}}dc_{1j}\int_{-\infty}^\infty \exp \left(-\frac{u_j^2}{2}-c_{1j}iu_j\right)\left(1-\left(\frac{3iu_j}{2}+\frac{c_{1j}u_j^2}{2}-\frac{iu_j^3}{3}\right)L^{-1/2}+O(L^{-1})\right)du_j \\
q_5 &\approx q\cdot I(c_{21},c_{12}) \\
q_6 &\approx 1-q_1-q\cdot I(c_{11},c_{12})
\end{aligned}
$$
where
\begin{align*}
I(\sigma_1,\sigma_2) = &-\frac{ q}{(2\pi i)^2}\int_{-\infty}^{\sigma_1}  d\sigma_1\int_{-\infty}^{\sigma_2}  d\sigma_2\int_{-\infty}^{\infty}   \int_{-\infty}^{\infty}\frac{u_1-u_2}{u_1-qu_2}\prod_{j=1}^2\exp \left(-\frac{u_j^2}{2}-\sigma_jiu_j\right)\\
&\quad \quad \quad \quad \times \left(1-\sum_{j=1}^2\left(\frac{3iu_j}{2}+\frac{\sigma_ju_j^2}{2}-\frac{iu_j^3}{3}\right)L^{-1/2}+O(L^{-1})\right)du_1du_2,
\end{align*}
where the $u_1$ and $u_2$ contours are deformed to avoid the pole at $u_1-qu_2$. 
\end{thm}

    

\section{Proof}

\subsection{Asymptotics of $q_1,q_2,q_3$}

Because the terms $q_1,q_2,q_3$ are simply the incomplete gamma function, an asymptotic analysis of the incomplete gamma function determines the asymptotics of $q_1,q_2$ and $q_3$. The full asymptotic expansion was determined in a recent paper \cite{ND19}. The result is that as $a\rightarrow \infty$, 



$$
Q\left(a, a+\tau a^{\frac{1}{2}}\right) \sim \frac{1}{2} \operatorname{erfc}\left(2^{-\frac{1}{2}} \tau\right)+\frac{1}{\sqrt{2 \pi a}} \exp \left(-\frac{\tau^{2}}{2}\right) \sum_{n=0}^{\infty} \frac{C_{n}(\tau)}{a^{n / 2}}.
$$
Somewhat awkwardly, the exponent of $1/2$ is in the second variable rather than the first. So if $L=a+\tau a^{1/2}$ then 
$$a={-}\tau\sqrt{ L +\tau^{2}/4}+ L+\tau^{2}/2= L {-} \tau\sqrt{L} + O(1).$$

Thus as $L\rightarrow\infty$,
$$
q_1 =Q((L+c_{11}L^{1/2} +O(1),L) \approx  \frac{1}{2} \operatorname{erfc}\left(- \frac{c_{11}}{2^{\frac{1}{2}}}\right)+\frac{1}{\sqrt{2 \pi a_1}} \exp \left(-\frac{c_{11}^{2}}{2}\right) \sum_{n=0}^{\infty} \frac{C_{n}\left(-c_{11}\right)}{a_1^{n / 2}},
$$
where $a_1=  L+\frac{c_{11}^2}{2}+c_{11}\sqrt{L+\frac{c_{i1}^2}{4}}   =L+c_{11}L^{1/2}+O(1)$.

\subsection{Asymptotics of $q_4,q_5,q_6$}

To analyze the asymptotics of $q_5$ and $q_6$, first denote 
\begin{multline*}
q(\sigma_1,\sigma_2)= \frac{ 1}{(2\pi i)^2}\int_{\tilde{C}_{1}}   \frac{dw_1}{w_1}  \int_{\tilde{C}_{2}}   \frac{dw_2}{w_2} \frac{w_1-w_2}{(w_1-qw_2)(1-w_1)(1-w_2)}\\
\times \prod_{j=1}^2\exp\left(LS\left(w_j\right)-L^{1 / 2}\sigma_j \log \left(1-w_j\right)\right)
\end{multline*}
where $S(w)=-w-\log(1-w)$. Then, $q_5=q\cdot q(c_{21},c_{12})$ and $q_6=1-q_1-q\cdot (c_{11},c_{12})$.
As observed in section 5 of \cite{Kua21}, $S(w)$ has a saddle point at $w=0$. Near $w=0$, we have the following Taylor expansions: 
$$S(u)=\sum_{j=2}^{\infty}\frac{u^j}{j},\quad\quad\log(1-u)=-\sum_{j=1}^{\infty}\frac{u^j}{j}.$$

Now change the contours to   vertical lines through $0$, let $w_i=-iL^{-1/2}u_i$
$$\begin{aligned}L S\left(w_j\right)-L^{1 / 2} \sigma_j \log \left(1-w_j\right)&=-\frac{u_j^2}{2}-\sigma_jiu_j-\left(\frac{\sigma_ju_j^2}{2}-\frac{iu_j^3}{3}\right)L^{-1/2}+O(L^{-1}).
\end{aligned}$$
To avoid the simple pole at $0$, we take partial derivatives with respect to $\sigma_j$, thus by \cite{Kua21} (The current arXiv version of \cite{Kua21} omits the comment about deforming the contours to avoid the pole, but it will be updated in the journal version)
$$\begin{aligned}
&\frac{\partial^2 q(\sigma_1,\sigma_2)}{\partial \sigma_1\partial \sigma_2}\\
&= \frac{ L}{(2\pi i)^2}\int \frac{\log(1-w_1)}{w_1}dw_1\int\frac{\log(1-w_2)}{w_2}dw_2\frac{w_1-w_2}{(w_1-qw_2)(1-w_1)(1-w_2)}\\
&\prod_{j=1}^2\exp\left(S\left(w_j\right)-L^{1 / 2}\sigma_j \log \left(1-w_j\right)\right)\\
& =\frac{ qL}{(2\pi i)^2}\int_{-\infty}^{\infty}   \int_{-\infty}^{\infty}\frac{u_1-u_2}{u_1-qu_2}\prod_{j=1}^2\frac{\log(1+iL^{-1/2}u_j)}{u_j(1+iL^{-1/2}u_j)}\\
&\exp \left(-\frac{u_j^2}{2}-\sigma_jiu_j-\left(\frac{\sigma_ju_j^2}{2}-\frac{iu_j^3}{3}\right)L^{-1/2}+O(L^{-1})\right)du_j\\
& =-\frac{ q}{(2\pi i)^2}\int_{-\infty}^{\infty}   \int_{-\infty}^{\infty}\frac{u_1-u_2}{u_1-qu_2}\prod_{j=1}^2\left(1-iL^{-1/2}u_j/2+O(L^{-1})\right)\left(1-iL^{-1/2}u_j+O(L^{-1})\right)\\
&\exp \left(-\frac{u_j^2}{2}-\sigma_jiu_j\right)\left(1-\left(\frac{\sigma_ju_j^2}{2}-\frac{iu_j^3}{3}\right)L^{-1/2}+O(L^{-1})\right)du_j\\
& =-\frac{ q}{(2\pi i)^2}\int_{-\infty}^{\infty}   \int_{-\infty}^{\infty}\frac{u_1-u_2}{u_1-qu_2}\prod_{j=1}^2\exp \left(-\frac{u_j^2}{2}-\sigma_jiu_j\right)\\
&\left(1-\sum_{j=1}^2\left(\frac{3iu_j}{2}+\frac{\sigma_ju_j^2}{2}-\frac{iu_j^3}{3}\right)L^{-1/2}+O(L^{-1})\right)du_1du_2.\\
\end{aligned}$$
Thus 
$$\begin{aligned}
  q(\sigma_1,\sigma_2)\approx& -\frac{ q}{(2\pi i)^2}\int_{-\infty}^{\sigma_1}  d\sigma_1\int_{-\infty}^{\sigma_2}  d\sigma_2\int_{-\infty}^{\infty}   \int_{-\infty}^{\infty}\frac{u_1-u_2}{u_1-qu_2}\prod_{j=1}^2\exp \left(-\frac{u_j^2}{2}-\sigma_jiu_j\right)\\
&\left(1-\sum_{j=1}^2\left(\frac{3iu_j}{2}+\frac{\sigma_ju_j^2}{2}-\frac{iu_j^3}{3}\right)L^{-1/2}+O(L^{-1})\right)du_1du_2.\\
\end{aligned}$$



Next analyze $q_4$. Deform contours in $q_4$ so that $0$ is outside of new contours:
$$\begin{aligned}
    q_4=&\frac{ 1}{(2\pi i)^2}\int_{C} \frac{d w_{1}}{w_{1}}  \int_{C} \frac{d w_{2}}{w_{2}} \frac{w_1-w_2}{w_1-qw_2} \prod_{j=1}^{2}\left[(1-w_j)^{-(y_j-x_1+1)}e^{-w_{j} t}\right]\\
    =&\frac{ 1}{(2\pi i)^2}\int_{C} \frac{d w_{1}}{w_{1}}  \int_{\tilde{C}_{2}} \frac{d w_{2}}{w_{2}} \frac{w_1-w_2}{w_1-qw_2} \prod_{j=1}^{2}\left[(1-w_j)^{-(y_j-x_1+1)}e^{-w_{j} t}\right]\\
    &+\frac{ 1}{2\pi i}\int_{C}  (1-w_1)^{-(y_1-x_1+1)}e^{-w_{1} t}\frac{d w_{1}}{w_{1}}\\
    =&q(c_{11},c_{12})+1+\frac{ 1}{2\pi i}\int_{\tilde{C}_{1}}  (1-w_1)^{-(y_1-x_1+1)}e^{-w_{1} t}\frac{d w_{1}}{w_{1}}\\
    &  +\frac{ 1}{2\pi i}\int_{\tilde{C}_{2}}  (1-w_2)^{-(y_2-x_1+1)}e^{-w_{2} t}\frac{d w_{2}}{w_{2}}.
\end{aligned}$$

Last, we take the limit as $L\xrightarrow[]{}\infty$ of the following integral:
$$\begin{aligned}
   & \frac{ 1}{2\pi i}\int_{\tilde{C}_{j}}  (1-w_j)^{-(y_j-x_1+1)}e^{-w_{j} t}\frac{d w_{j}}{w_{j}(1-w_j)}\\
   &=\frac{ 1}{2\pi i}\int_{\tilde{C}_{j}}   \exp \left(L S\left(w_j\right)-L^{1 / 2} c_{1j} \log \left(1-w_j\right)\right)\frac{d w_{j}}{w_{j}}\\
   &=\frac{ -L^{1/2}}{2\pi i}\int_{-\infty}^{c_{1j}}dc_{1j}\int_{\tilde{C}_{j}} \frac{\log(1-w_j)}{w_{j}(1-w_j)}  \exp \left(L S\left(w_j\right)-L^{1 / 2} c_{1j} \log \left(1-w_j\right)\right)d w_{j}\\
    &\approx\frac{1 }{2\pi }\int_{-\infty}^{c_{1j}}dc_{1j}\int_{-\infty}^\infty \exp \left(-\frac{u_j^2}{2}-c_{1j}iu_j\right)\left(1-\left(\frac{3iu_j}{2}+\frac{c_{1j}u_j^2}{2}-\frac{iu_j^3}{3}\right)L^{-1/2}+O(L^{-1})\right)du_j. \\
\end{aligned} $$


\end{document}